\documentclass[11pt]{amsart}
\usepackage{pictex}
\usepackage{amssymb,amsmath,amsgen,amsxtra,amsfonts,a4wide} 
\usepackage{dsfont,times}

\begin{document}
%
%
%
\theoremstyle{definition}
\newtheorem{Definition}{Definition}[section]
\newtheorem{Construction}{Construction}[section]
\newtheorem{Example}[Definition]{Example}
\newtheorem{Examples}[Definition]{Examples}
\newtheorem{Remark}[Definition]{Remark}
\newtheorem{Remarks}[Definition]{Remarks}
\newtheorem{Caution}[Definition]{Caution}
\newtheorem{Conjecture}[Definition]{Conjecture}
\newtheorem{Question}[Definition]{Question}
\newtheorem{Questions}[Definition]{Questions}
\theoremstyle{plain}
\newtheorem{Theorem}[Definition]{Theorem}
\newtheorem*{Theoremx}{Theorem}
\newtheorem{Proposition}[Definition]{Proposition}
\newtheorem{Lemma}[Definition]{Lemma}
\newtheorem{Corollary}[Definition]{Corollary}
\newtheorem{Fact}[Definition]{Fact}
\newtheorem{Facts}[Definition]{Facts}
\newtheoremstyle{voiditstyle}{3pt}{3pt}{\itshape}{\parindent}%
{\bfseries}{.}{ }{\thmnote{#3}}%
\theoremstyle{voiditstyle}
\newtheorem*{VoidItalic}{}
\newtheoremstyle{voidromstyle}{3pt}{3pt}{\rm}{\parindent}%
{\bfseries}{.}{ }{\thmnote{#3}}%
\theoremstyle{voidromstyle}
\newtheorem*{VoidRoman}{}

%
\newcommand{\prf}{\par\noindent{\sc Proof.}\quad}
\newcommand{\blowup}{\rule[-3mm]{0mm}{0mm}}
\newcommand{\cal}{\mathcal}
\newcommand{\Aff}{{\mathds{A}}}
\newcommand{\BB}{{\mathds{B}}}
\newcommand{\CC}{{\mathds{C}}}
\newcommand{\FF}{{\mathds{F}}}
\newcommand{\GG}{{\mathds{G}}}
\newcommand{\HH}{{\mathds{H}}}
\newcommand{\NN}{{\mathds{N}}}
\newcommand{\ZZ}{{\mathds{Z}}}
\newcommand{\PP}{{\mathds{P}}}
\newcommand{\QQ}{{\mathds{Q}}}
\newcommand{\RR}{{\mathds{R}}}
\newcommand{\Liea}{{\mathfrak a}}
\newcommand{\Lieb}{{\mathfrak b}}
\newcommand{\Lieg}{{\mathfrak g}}
\newcommand{\Liem}{{\mathfrak m}}
\newcommand{\ideala}{{\mathfrak a}}
\newcommand{\idealb}{{\mathfrak b}}
\newcommand{\idealg}{{\mathfrak g}}
\newcommand{\idealm}{{\mathfrak m}}
\newcommand{\idealI}{{\cal I}}
\newcommand{\lin}{\sim}
\newcommand{\num}{\equiv}
\newcommand{\dual}{\ast}
\newcommand{\iso}{\cong}
\newcommand{\homeo}{\approx}
\newcommand{\mm}{{\mathfrak m}}
\newcommand{\pp}{{\mathfrak p}}
\newcommand{\qq}{{\mathfrak q}}
\newcommand{\rr}{{\mathfrak r}}
\newcommand{\pP}{{\mathfrak P}}
\newcommand{\qQ}{{\mathfrak Q}}
\newcommand{\rR}{{\mathfrak R}}
%
%
\newcommand{\dq}{{``}}
\newcommand{\OO}{{\cal O}}
\newcommand{\into}{{\hookrightarrow}}
\newcommand{\onto}{{\twoheadrightarrow}}
\newcommand{\Spec}{{\rm Spec}\:}
\newcommand{\Pic}{{\rm Pic }}
\newcommand{\Br}{{\rm Br}}
\newcommand{\NS}{{\rm NS}}
\newcommand{\chit}{\chi_{\rm top}}
\newcommand{\KanDiv}{{\cal K}}
\newcommand{\perdef}{{\stackrel{{\rm def}}{=}}}
\newcommand{\Cycl}[1]{{\ZZ/{#1}\ZZ}}
\newcommand{\Sym}{{\mathfrak S}}
\newcommand{\Alt}{{\mathfrak A}}
\newcommand{\ab}{{\rm ab}}
\newcommand{\Aut}{{\rm Aut}}
\newcommand{\Hom}{{\rm Hom}}
\newcommand{\Supp}{{\rm Supp}}
\newcommand{\ord}{{\rm ord}}
\newcommand{\Alb}{{\rm Alb}}
\newcommand{\Jac}{{\rm Jac}}
\newcommand{\piet}{{\pi_1^{\rm \acute{e}t}}}
\newcommand{\Het}[1]{{H_{\rm \acute{e}t}^{{#1}}}}
\newcommand{\Hcris}[1]{{H_{\rm cris}^{{#1}}}}
\newcommand{\HdR}[1]{{H_{\rm dR}^{{#1}}}}
\newcommand{\hdR}[1]{{h_{\rm dR}^{{#1}}}}
\newcommand{\defin}[1]{{\bf #1}}
\title[Uniruled Surfaces of General Type]{Uniruled Surfaces of General Type}
\author{Christian Liedtke}
\address{Mathematisches Institut, Heinrich-Heine-Universit\"at, 40225
  D\"usseldorf, Germany}
\email{liedtke@math.uni-duesseldorf.de}
\date{November 6, 2006}

\begin{abstract}
  We give a systematic construction of uniruled surfaces 
  in positive characteristic.
  Using this construction,
  we find surfaces of general type with non-trivial global
  vector fields, surfaces with arbitrarily non-reduced Picard schemes
  as well as surfaces with inseparable canonical maps.
  In particular, we show that some previously known
  pathologies are not sporadic but exist in abundance.
\end{abstract}
\setcounter{tocdepth}{1}
\maketitle
\tableofcontents
\section*{Introduction}

In the beginning of the 20th century, 
the Italian school established a coarse classification
of complex surfaces of special type, the so-called
Enriques classification.
Since then, it has been clarified and refined by
Kodaira, \v{S}afarevi\v{c} and many others.
In a series of three papers in the 1970s, 
Bombieri and Mumford \cite{bm3} 
extended this classification 
to positive characteristic.
A major new feature is the existence of fibrations
over curves with singular geometric generic fibres.
For surfaces of special type, only quasi-elliptic
surfaces, 
which exist in characteristic $2$ and $3$ only,
have to be considered.
By definition, a quasi-elliptic surface is a surface that
admits a fibration over a curve such that the 
geometric generic fibre
is a singular rational curve of arithmetic genus $1$.
In particular, a quasi-elliptic surface is (inseparably)
uniruled.

But what about the classification of complex surfaces of
general type?
A full classification seems out of range at the moment.
However, one would like to have bounds on their
invariants, to understand the behaviour of their 
pluricanonical systems, and to classify surfaces with low invariants.
Famous results in this direction are the Bogomolov--Miyaoka--Yau
inequality $c_1^2\leq 9\chi$ and Bombieri's analysis of the
pluricanonical systems, just to name two.
The situation in positive characteristic is more complicated: 
many equalities and inequalities, which easily follow 
from Hodge theory, do not hold simply because the 
Fr\"olicher spectral sequence from Hodge to de~Rham cohomology
does not degenerate at $E_1$-level in general.
The Bogomolov--Miyaoka--Yau inequality is known to 
fail \cite[Section 3.4.1]{sz}.
On the other hand, Ekedahl \cite{ek2} and Shepherd-Barron \cite{sb}
extended Bombieri's results on pluricanonical maps to positive
characteristic.
Again, the exceptions occur mainly in small characteristics
and (inseparably) uniruled surfaces are responsable 
for many unexpected phenomena.
Thus, in order to understand failures of classical theorems 
about surfaces in positive characteristic, it is indispensable to 
understand uniruled surfaces.
And as a first step, one should have at least a large
supply of examples to study.
\medskip

In this paper, we present a systematic construction of
surfaces that are uniruled or
birationally dominated by Abelian surfaces.
This construction is inspired by the structure result of 
quasi-hyperelliptic surfaces \cite[Theorem 1]{bm3}
and Schr\"oer's construction \cite{sch} of unirational 
K3 surfaces.

To obtain our surfaces, we need two curves $C,F$ 
and rational $p$-closed vector fields $\delta_C, \delta_F$ 
on them.
Resolving the singularities of the quotient of 
$C\times F$ by $\delta_C+\delta_F$, we obtain a surface $X$
which is birationally dominated by $C\times F$.
In particular, if at least one of the curves is rational,
the surface $X$ is uniruled.
If both curves are rational then $X$ is even unirational.
The surface $X$ comes with fibrations
over $C^{(-1)}$ and $F^{(-1)}$, and
both fibrations are usually not generically smooth.

We will restrict us mostly to characteristic $2$.
This is mainly because we expect to find more pathologies 
in this characteristic and partly to simplify our
exposition.
For example, the restriction to characteristic $2$ allows
an analysis of the slope spectral sequence and the
spectral sequence from Hodge to de~Rham cohomology.
However, before stating these results 
we first discuss the examples we have found.
\medskip

For curves of general type, i.e. of genus at least $2$, the
canonical map is always a separable morphism, which is
either of degree $2$ onto $\PP^1$ or defines an embedding.
Already in characteristic zero,
the situation for surfaces is more complicated and
we refer to \cite[Section VII.7]{bhpv} for an introduction 
as well as references.
Here we show that the canonical map can become inseparable
and that this is not a sporadic phenomenon:

\begin{VoidItalic}[Theorem \ref{canonicalinsep}]
    In characteristic $2$ there exist
    unbounded families of unirational
    surfaces of general type whose
    canonical maps are inseparable morphisms onto rational
    surfaces.
\end{VoidItalic}

We note that it follows from Shepherd-Barron's results 
\cite[Theorem 27]{sb}
that $|3K_X|$ of a surface of general type
defines a birational morphism provided that $c_1^2$
and $\chi$ are sufficiently large.
Moreover, if $X$ does not possess a pencil of curves of
arithmetic genus $2$ then
already $|2K_X|$ defines a birational morphism if
$c_1^2$ and $\chi$ are sufficiently large.
However, we can arrange the surfaces of Theorem \ref{canonicalinsep}
not to possess pencils of curves of small arithmetic genus.
Hence the inseparability of the canonical map is not related
to the existence of special fibrations of low genus.\medskip

It is already known for some time that
the Bogomolov--Miyaoka--Yau inequality $c_1^2\leq9\chi$ may
fail in positive characteristic \cite[Section 3.4.1]{sz}.
Here we present a surface violating this inequality with
$\chi=1$, which is the smallest value possible for $\chi$
in characteristic zero.

\begin{VoidItalic}[Theorem \ref{bmy}]
  In characteristic $2$ there exist
  surfaces of general type with $\chi=1$ and $c_1^2=14$.
\end{VoidItalic}
\medskip

The number of isolated $(-2)$-curves on a minimal surface of
general type is bounded above by $\frac{1}{9}(3c_2-c_1^2)$ 
by a theorem of Miyaoka.
Also this is known to fail in positive characteristic
and Shepherd-Barron \cite[Theorem 4.1]{sbfol} has shown that
if the number of $(-2)$-curves exceeds 
$c_1^2+\frac{1}{2}c_2$ then the surface is uniruled.
However, usually there is a gap between these two bounds
and we show that this gap is populated by uniruled as well
at non-uniruled surfaces.
In particular, Shepherd-Barron's bound is not sharp:

\begin{VoidItalic}[Theorem \ref{minustwo}]
  There exist minimal surfaces of general type in characteristic $2$
  violating Miyaoka's bound on $(-2)$-curves that do not reach 
  Shepherd-Barron's bound.
  There exist uniruled as well as non-uniruled
  such surfaces.
\end{VoidItalic}
\medskip

Since a group scheme over a field of positive characteristic
may be non-reduced one has to distinguish between the Picard
variety and the Picard scheme of a variety.
Examples of surfaces with non-reduced Picard schemes 
fields are known, e.g. \cite{ig2}.
However, one could ask whether there are bounds on the
non-reducedness, e.g. one could ask whether the dimension
$h^1(\OO_X)$ of the tangent space to the Picard scheme 
$\Pic(X)$ is bounded, say, in terms of the dimension
$\frac{1}{2}b_1(X)$ of the Picard scheme.
This is not the case:

\begin{VoidItalic}[Theorem \ref{picard}]
  Given an integer $q\geq2$, there exists a family 
  $\{X_i\}_{i\in\NN}$ of
  uniruled surfaces of general type in
  characteristic $2$
  all having the same Picard variety 
  of dimension $q$ such that
  $$
   h^{01}(X_i)\,=\,h^1(\OO_{X_i})\,\to\,\infty
     \,\mbox{ as }\,i\to\infty
  $$
  Thus, the Picard scheme can get arbitrarily non-reduced,
  even when fixing the Picard variety.
\end{VoidItalic}
\medskip

It follows from Hodge theory that all global $1$-forms on a 
complex projective manifold are pull-backs of global $1$-forms from
its Albanese variety via the Albanese map.
On the other hand, Igusa \cite{ig2} gave an example of a surface with 
$h^0(\Omega_X^1)$ strictly larger than 
$\frac{1}{2}b_1(X)=\dim\Alb(X)=h^0(\Omega_{\Alb(X)}^1)$ and 

\begin{VoidItalic}[Theorem \ref{albanese}]
  Given an integer $q\geq2$, there exists a family 
  $\{X_i\}_{i\in\NN}$ of
  uniruled surfaces of general type in
  characteristic $2$
  all having the same Albanese variety 
  of dimension $q$ such that
  $$
   h^{10}(X_i)\,=\,h^0(\Omega_{X_i}^1)\,\to\,\infty
   \,\mbox{ as }\,i\to\infty\,.
  $$
\end{VoidItalic}

Combining Theorem \ref{picard} and Theorem \ref{albanese} we can
even produce families $\{X_i\}_{i\in\NN}$ of uniruled surfaces 
of general type with fixed Albanese variety and
$h^{10}(X_i)\,-\,h^{01}(X_i)$ tending to infinity as $i$ tends
to infinity.
Thus, we can also violate the Hodge symmetry \dq$h^{10}\,=\,h^{01}$\dq\ 
as much as we want - even when fixing the first Betti number.
\medskip

The automorphism group of a surface
of general type over the complex numbers is finite.
Hence its Lie algebra, which can be identified with
the space of global holomorphic vector fields, is
trivial.
In particular, a surface of general type over the
complex numbers does not possess non-trivial global
vector fields.
However, surfaces of general type with non-trivial
global vector fields in positive 
characteristic have been constructed by 
Lang \cite{lang}, Shepherd-Barron \cite{sbfol} 
and others.
Again, we are able to obtain this phenomenon in
unbounded families.

\begin{VoidItalic}[Theorem \ref{unboundedvectorfields}]
  In characteristic $2$, there exist unbounded families of 
  surfaces of general
  such that each member
  of this family possesses non-trivial global vector fields.
  Moreover, we can find such families in which every member is
  uniruled, resp. not uniruled.
\end{VoidItalic}
\medskip

In characteristic $2$, our construction has the nice 
feature that it is possible to compute invariants that are
usually difficult to determine.
Since we can easily produce surfaces that do not lift
to characteristic zero, even over a ramified extension of
the Witt ring, these computations may be an interesting
testing ground for general conjectures about surfaces.

More precisely, we determine the Betti and 
Hodge numbers and analyse the spectral sequences related to
Hodge and de~Rham--Witt cohomology.
It turns out that the spectral sequence
from Hodge to de~Rham cohomology may or may not degenerate
at $E_1$-level (we present examples for both cases),
whereas  the slope spectral sequence from de~Rham--Witt to 
crystalline cohomology usually does not degenerate
at $E_1$-level.
Also, we verify the Artin--Tate conjecture for our surfaces and
show that they possess an Artin invariant, analogous to
the one defined by Artin \cite{ar}
for supersingular K3 surfaces. 

\begin{VoidRoman}[Acknowledgements]
  I thank Stefan~Schr\"oer for many long and stimulating
  discussions.
  Also, I thank Frans~Oort and Kay~R\"ulling for pointing out
  inaccuracies.
\end{VoidRoman}

\section{Rational vector fields on curves}
\label{curvesection}

This section deals with rational vector fields on smooth
curves over fields of characteristic $p>0$.
We need a supply of such vector fields for our constructions.
But since a smooth curve of genus $g\geq2$ does not possess any
non-trivial vector fields and the existing
vector fields on curves of genus $0$ and $1$ are usually
not interesting for us, we have to 
work with rational vector fields from the very 
beginning.

For a vector field $\delta_C$ on a curve $C$ we denote
by $(\delta_C)$ its divisor,
by $(\delta_C)_0$ its divisor of zeros and by
$(\delta_C)_\infty$ its divisor of poles.
Thus, $(\delta_C)=(\delta_C)_0-(\delta_C)_\infty$.

We recall that a rational vector field $\delta$ is called 
{\em $p$-closed}, if $\delta^{[p]}=f\cdot\delta$ for some
rational function $f$ on $C$.
If $f=0$ the vector field is called {\em additive}, whereas
it is called {\em multiplicative} if $f=1$.

\subsection*{Rational curves}
Let $x$ be a coordinate on $\PP^1$.
We consider the following rational vector fields,
which are easily seen to be additive 
in characteristic $2$.
\begin{eqnarray}
\label{vectorone}
\delta_1 &:=& (x^{-4}+x^{-2})  D_x \\
\label{vectortwo}
\delta_2 &:=& (x^{-2} + x^4)  D_x
\end{eqnarray}
The zeros of both vector fields are of order $2$.
The vector field $\delta_1$ has a pole of order $4$ at $x=0$, whereas
$\delta_2$ has poles of order $2$ at $x=0$ and $x=\infty$.

More generally, we choose pairwise distinct elements
$a_1,...,a_n$ and $b_1,...,b_n$ of the ground field $k$.
Then the rational vector field
\begin{eqnarray}
  \label{singrational}
  \delta_n' &:=& \prod_{i=1}^n (x-a_i)^2\, (x-b_i)^{-2} \,D_x 
\end{eqnarray}
has only zeros and poles of order $2$ and is additive in
characteristic $2$.
More precisely, its poles lie at $x=b_i$ and its $n+1$ zeros
are given by $x=a_i$ and $x=\infty$.

\subsection*{Elliptic curves}
The Deuring normal form of an elliptic
curve in characteristic $2$ is given by
$$  
  y^2\,+\alpha xy\,+\,\,y\,=\,x^3\,,
$$
where $\alpha\in k$ satisfies $\alpha^3\neq1$, cf.
\cite[Appendix A]{silv}.
We denote by $E_\alpha$ the closure of this affine curve
in $\PP^2$, which is a smooth elliptic curve
with $j$-invariant $\alpha^{12}/(\alpha^3-1)$.
There exists a $p$-closed 
regular vector field on $E_\alpha$ of additive type, 
i.e. the curve is supersingular, if and
only if $\alpha=0$.

The rational vector field 
\begin{eqnarray}
  \label{singelliptic}
\delta_{\alpha,a,b} &:=&
   ( a \, + \, b\cdot x ) \,\cdot\, 
   ( (1+\alpha x)\, D_x + (\alpha y+x^2)\, D_y ) \,.
\end{eqnarray}
on $\PP^2$ descends to a rational vector field
on $E_\alpha$.
This rational vector field is additive iff $a\alpha+b=0$ and
multiplicative iff $a\alpha+b=1$.
In particular, if $\alpha\neq0$, i.e. if $E_\alpha$ is not
supersingular, the rational vector field
$\delta_{\alpha,1,\alpha}$
on $E_\alpha$ is additive and has one zero and one
pole of order $2$.

We would like to mention how we found these vector fields:
If we consider the affine curve in Deuring normal form
as lying inside $\Aff^1\times\Aff^1$ then its projective
closure $F_\alpha$ in $\PP^1\times\PP^1$ 
is a singular elliptic curve with a cusp at infinity.
The vector space $H^0(F_\alpha,\Theta_{F_\alpha})$ is
$2$-dimensional, where $\Theta_{F_\alpha}$ denotes the
dual of the sheaf of K\"ahler differentials.
This space is explicitly described by (\ref{singelliptic})
and a rational vector field extends to a regular vector
field on its normalisation $E_\alpha$ iff $a=0$.

\subsection*{Hyperelliptic curves in characteristic 
$\mathbf2$}
A smooth curve $C$ of genus $g\geq2$ 
is called {\em hyperelliptic} if it 
admits a separable morphism $\varphi:C\to\PP^1$ of
degree $2$.
By \cite[Section 1]{bh} we may assume that $\varphi$ is
branched over $g+1$ distinct points $P_1,...,P_{g+1}$.
In particular, all higher ramification groups $G_i$ vanish
already for $i\geq2$, which is the lowest value possible 
in presence of wild ramification, cf. 
\cite[Chapitre IV]{se}.
We let $x$ be a parameter on $\PP^1$. 
By \cite[Proposition 1.5]{bh}, there exists a polynomial
$g(x)$ of degree $g+1$ that does not vanish in any of the
$P_i$ such that $C$ is given over 
$\Aff^1:=\PP^1-\{x=\infty\}$ by the equation
$$
  z^2\,+\,f(x)\,z\,+\,f(x)\,g(x)\,=\,0,
  \mbox{ where }
  f(x)\,:=\,\prod_{i=1}^{g+1}(x-\alpha_i)\,,
$$
where the $\alpha_i$ correspond to the $P_i$.
Straight forward local calculations give the divisor
\begin{eqnarray}
\label{singhyperelliptic}
   (\partial/\partial x) &=& 
   2\,P_\infty' \,+\,2\,P_\infty'' \,+\, 
   \sum_{i=1}^{g+1} (-2)\,P_i'  \,.
\end{eqnarray}
Here, $P_i'$ is the unique point of $C$ lying above $P_i$,
and $P_\infty'$, $P_\infty''$ lie above $x=\infty$.

\subsection*{Artin--Schreier extensions of $\mathbf\PP^1$}
We consider the Artin--Schreier extension
$$
z^p\,-\,z\,=\,x^{hp-1}
$$
of the affine line in characteristic $p$.
Its projective closure is a curve $C$ of genus
$g=1-p+\frac{1}{2}p(p-1)h$ together
with a morphism $C\to\PP^1$, which is wildly ramified 
at infinity.
Pulling back $x$ we obtain
\begin{equation}
 \label{singartin}
 (\partial/\partial x) \,=\, p(h(p-1)-2)\,P_\infty \,.
\end{equation}
By the Deuring-\v{S}afarevi\v{c}-formula 
\cite[Corollary 1.8]{cr},
the $p$-rank of such a curve is equal to zero.
In particular, a line bundle $\cal L$ on $C$
with ${\cal L}^{\otimes p}\iso\OO_C$ 
is trivial.

\section{Singular vector fields on surfaces}

In this section we consider rational $p$-closed
vector fields on surfaces with small multiplicity.
In characteristic $p=2$, this allows us to determine the
the singularities of the quotient of a surface by these
vector fields.

First, we recall some well-known facts from \cite[Section \S1]{rs}.
On a smooth surface $S$, we can write a rational vector field
around a point $P$ in local coordinates as
$$\delta\,=\,
h(x,y) \,\left(\, f(x,y)\,\frac{\partial}{\partial x}\,+\, 
g(x,y)\,\frac{\partial}{\partial y} \,\right)\,
$$
where $h$ is a rational function and
where $f$ and $g$ are regular functions around the point 
$P$ such that the ideal $\idealI:=(f,g)$ generated by $f$ and $g$ 
has height at least $2$.

If $\idealI$ is not the unit ideal, the vector field is said to
have an {\em isolated singularity} at $P$.
The {\em multiplicity} of $\delta$ in $P$ is the dimension
of the $k$-vector space $\OO_{P}/\idealI$.
In case the rational vector field has no isolated singularities
the vector field is said to have 
{\em only divisorial singularities}.
Locally around $P$, the function $h$ 
defines a divisor and all these functions at all
points of $S$ define a divisor, the
{\em divisor} $(\delta)$ of the rational vector field.

If $\delta$ is a $p$-closed vector field on $S$ then
we can form the quotient $S/\delta$, which is a normal surface.
Its isolated singularities lie below those points of $S$ where
the vector field $\delta$ has an isolated singularity.
If we blow up $S$ at an isolated singularity of $\delta$,
then we obtain an induced vector field on the blow-up.
One would like to find a finite sequence of blow-ups such that
the induced vector field $\tilde{\delta}$ on the blow-up 
$\tilde{S}$ is a rational vector
field with only divisorial singularities.
This would yield a diagram
$$\begin{array}{ccc}
   S & \leftarrow & \tilde{S} \\
   \downarrow & &\downarrow\\
   S/\delta &\leftarrow & \tilde{S}/\tilde{\delta}
 \end{array}
$$
where $\tilde{S}/\tilde{\delta}$ resolves the singularities of
$S/\delta$.
In general, this is not possible.
However, there is the following remarkable result
from \cite[Proposition 2.6]{hir}.
 
\begin{Proposition}[Hirokado]
  \label{hirokado}
  The singularities of a $p$-closed vector field on a surface
  in characteristic $p=2$
  can be resolved by a sequence of blow-ups.
\end{Proposition}

We recall that a singularity on $X$ 
is called {\em rational} (resp. {\em elliptic}) 
if $R^1 f_\ast\OO_{\tilde{X}}$ 
is a zero-dimensional (resp. one-dimensional) 
vector space for one, and hence every,
resolution of singularities $f:\tilde{X}\to X$.
In our examples, we will need the following 
singularities and their dual resolution graphs:
every vertex represents a rational curve, which has
self-intersection number $-2$ unless
labelled differently.

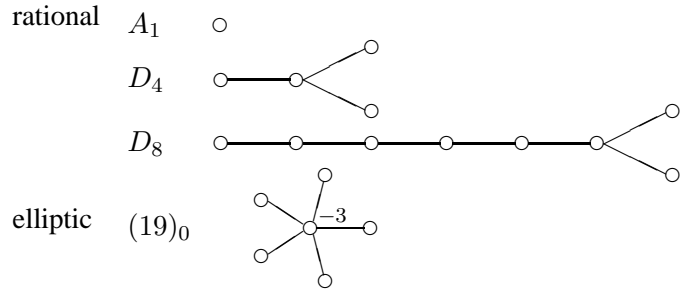
\begin{center}
\begin{tabular}{ll}
\blowup rational & \beginpicture
\setcoordinatesystem units <10mm,5mm>
\unitlength1pt
\put {$A_1$} at -1 0
\put {\circle{5}} [Bl] at 0 0
\endpicture \\
\blowup & \beginpicture
\setcoordinatesystem units <10mm,5mm>
\unitlength1pt
\put {$D_4$} at -1 0
\put {\circle{5}} [Bl] at 0 0
\put {\circle{5}} [Bl] at 1 0
\put {\circle{5}} [Bl] at 2 0.85
\put {\circle{5}} [Bl] at 2 -0.85
\put {\line(1,0){23.5}} [Bl] at 0.09 0
\put {\line(2,1){23.5}} [Bl] at 1.09 0
\put {\line(2,-1){23.5}} [Bl] at 1.09 0
\endpicture
 \\
\blowup & \beginpicture
\setcoordinatesystem units <10mm,5mm>
\unitlength1pt
\put {$D_8$} at -1 0
\put {\circle{5}} [Bl] at 0 0
\put {\circle{5}} [Bl] at 1 0
\put {\circle{5}} [Bl] at 2 0
\put {\circle{5}} [Bl] at 3 0
\put {\circle{5}} [Bl] at 4 0
\put {\circle{5}} [Bl] at 5 0
\put {\circle{5}} [Bl] at 6 0.85
\put {\circle{5}} [Bl] at 6 -0.85
\put {\line(1,0){23.5}} [Bl] at 0.09 0
\put {\line(1,0){23.5}} [Bl] at 1.09 0
\put {\line(1,0){23.5}} [Bl] at 2.09 0
\put {\line(-1,0){23.5}} [Bl] at 3.91 0
\put {\line(1,0){23.5}} [Bl] at 4.09 0
\put {\line(2,1){23.5}} [Bl] at 5.09 0
\put {\line(2,-1){23.5}} [Bl] at 5.09 0
\endpicture \\
\blowup elliptic & \beginpicture
\setcoordinatesystem units <10mm,5mm>
\unitlength1pt
\put {$(19)_0$} at -1 0
\put {\circle{5}} [Bl] at 0.35 0.75
\put {\circle{5}} [Bl] at 0.35 -0.75
\put {\circle{5}} [Bl] at 1 0
\put {\circle{5}} [Bl] at 1.8 0
\put {\circle{5}} [Bl] at 1.2 1.4
\put {\circle{5}} [Bl] at 1.2 -1.4

\put {$\scriptstyle -3$} at 1.3 0.3

\put {\line(1,0){18}} [Bl] at 1.09 0
\put {\line(1,4){3.8}} [Bl] at 1.05 0.17
\put {\line(1,-4){3.8}} [Bl] at 1.05 -0.17
\put {\line(3,-2){13}} [Bl] at 0.45 0.7
\put {\line(3,2){13}} [Bl] at 0.45 -0.7

\endpicture \\
\blowup &\\
\end{tabular}
\end{center}

Let us now assume that the vector field locally takes the form
\begin{equation}
  \label{vectorform}
  \delta \,=\, f(x)\,\frac{\partial}{\partial x}\,+\, 
  g(y)\,\frac{\partial}{\partial y}\,,
\end{equation}
where $f$ and $g$ are rational functions.
From the preceding discussion it is clear that 

\begin{Remark}
  \label{singularremark}
  The isolated singularities of the vector field
  (\ref{vectorform}) are those points where $f$ and
  $g$ both have a pole or where they both have a zero.  
\end{Remark}

We now determine the type of singularity the quotient acquires
if the multiplicity is small.

\begin{Proposition}
  \label{singularprop}
  Let $\delta$ be a $p$-closed rational vector field in characteristic
  $2$ as in formula (\ref{vectorform}).
  
  If $\delta$ has an isolated singularity at the origin $x=y=0$ then the
  quotient acquires
  \begin{enumerate}
  \item a rational singularity of type $A_1$ if $|\ord_x f|=|\ord_y g|=1$,
  \item a rational singularity of type $D_4$ if $|\ord_x f|=|\ord_y g|=2$,
  \item a rational singularity of type $D_8$ if $|\ord_x f|=4$ and $|\ord_y g|=2$,
  \item an elliptic singularity of type $(19)_0$ if $|\ord_x f|=|\ord_y g|=4$,
  \end{enumerate}
  at the point lying below the origin.
\end{Proposition}

\prf
In the first case, the vector field has multiplicity $1$ and it follows
from \cite[Proposition 2.2]{hir} and \cite[Corollary 2.5]{hir} 
that the quotient is a Du~Val singularity of type $A_1$.

Now suppose that $\ord_x f=\ord_y g=2$.
Then there exist regular functions 
$\epsilon(x)=\epsilon_0+\epsilon_1 x+...$
and $\eta(y)=\eta_0+\eta_1y+...$ with $\epsilon_0\neq0$
and $\eta_0\neq0$ such that
$$
\delta \,=\, x^2\epsilon(x)\,\frac{\partial}{\partial x}\,+\,
y^2\,\eta(y)\,\frac{\partial}{\partial y}\,.
$$
On the blow-up with coordinates $s$ and $y$, where 
$x=sy$, the induced rational vector field is
$$
 \tilde{\delta}\,=\,
 y\,\left(\, ( s^2\epsilon(sy) - s\eta(y))\,\frac{\partial}{\partial s}
   \,+\,y\eta(y)\,\frac{\partial}{\partial y}\,   \right) \,,
$$
where $y=0$ is the local equation of the exceptional divisor $E$ of
the blow-up.

The vector field $\tilde{\delta}$ has three isolated singularities
of multiplicity $1$ on $E$ at 
$s=0$, $s=\infty$ and $s=\eta_0/\epsilon_0$.
By Hirokado's result quoted above, these correspond to 
$A_1$-singularities and blowing up these three isolated points
the induced vector field $\tilde{\delta}$ 
on this blow-up $\tilde{S}$ has only divisorial singularities.
This yields the resolution graph of the singularity which
looks like $D_4$.

The exceptional divisors of the resolution of singularities
$\tilde{S}/\tilde{\delta}\to S/\delta$ are dominated by
the exceptional divisors of the blow-up $\tilde{S}\to S$.
In particular, all exceptional divisors of the resolution
of singularities are rational curves.
Using \cite[\S1, Proposition 1]{rs}, it is easy to see that
all self-intersection numbers are equal to $(-2)$.
Since the intersection of these curves with the canonical
class turns out to be zero, these curves are in fact smooth
rational curves.
This shows that the quotient by $\delta$ is indeed a singularity
of type $D_4$.

We leave the tedious calculations of the remaining cases 
to the reader.
\qed

\section{Uniruled surfaces}

We now present our construction for uniruled surfaces.
It is inspired by the classification of 
quasi-hyperelliptic surfaces by Bombieri and Mumford
\cite[Section 2]{bm3},
as well as the construction
of K3 surfaces via the self-product of two
cuspidal rational curves
by Schr\"oer \cite{sch}.\medskip

We consider the following data $(C,\,F,\,\delta\,:=\,\delta_C+\delta_F)$
\begin{enumerate}
  \item two smooth curves $C$ and $F$,
  \item a $p$-closed rational vector field $\delta_C$ on $C$ and
  \item a $p$-closed rational vector field $\delta_F$ on $F$, where
  \item $\delta_C$ and $\delta_F$ are either both additive or both 
    multiplicative.
\end{enumerate}
We define $S:=C\times F$ and denote by $\delta$ the
vector field $\delta_C+\delta_F$ on $S$.
The vector field $\delta$ is locally of the form
(\ref{vectorform}) and so Remark \ref{singularremark}
as well as Proposition \ref{singularprop} apply.

A theorem of Jacobson \cite[Section II.7.2]{dg} states
that 
$$
 (\delta_C\,+\delta_F)^{[p]}\,=\, \delta_C^{[p]} \,+\, 
 \Lambda_p (\delta_C,\delta_F) \,+\,\delta_F^{[p]}\,,
$$
where $\Lambda_p(-,-)$ is a universal expression in terms
of iterated Lie brackets.
Considered as vector fields on $S$, the Lie bracket 
$[\delta_C,\delta_F]$ is equal to zero and so
also $\Lambda_p(\delta_C,\delta_F)$ is zero.
Hence, if $\delta_C$ and $\delta_F$ are both additive 
(resp. multiplicative) the same is true 
for $\delta$.

\begin{Definition}
  \label{datadef}
  Given $(F,C,\delta)$ we let $S:=C\times F$.
  We will say that $X$ is a
  {\em surface constructed from data $(C,F,\delta)$} 
  if there exists a sequence of blow-ups $\tilde{S}\to S$ such that
  the induced vector field $\tilde{\delta}$ on $\tilde{S}$ has only
  divisorial singularities such that $X=\tilde{S}/\tilde{\delta}$ .
  If $F$ is a rational curve we will refer to $X$ as
  a {\em uniruled surface constructed from data $(C,F,\delta)$}.
\end{Definition}

Indeed, if $F$ is a rational curve then $X$ is dominated by 
a blow-up of a ruled surface and thus $X$ is (inseparably)
uniruled.
In case $C$ and $F$ are both rational, the surface $X$
is (inseparably) unirational.

Given a surface $X$ constructed from data $(C,F,\delta)$,
we have a commutative diagram
$$\begin{array}{ccc}
   S & \leftarrow & \tilde{S} \\
   \downarrow & &\downarrow\\
   S/\delta &\leftarrow & X
 \end{array}
$$
where $X$ is a smooth surface since $\tilde{\delta}$ has only
divisorial singularities.
More precisely, $X$ resolves the singularities of $S/\delta$.
The map $\tilde{S}\to X$ is a finite and purely inseparable morphism
of degree $p$ and height $1$, which makes the computation of the 
invariants of $X$ quite easy as we will see.
We note that by Hirokado's result 
(see Proposition \ref{hirokado} above),
the assumption on the existence of a suitable blow-up
$\tilde{S}\to S$ in Defintion \ref{datadef}
is automatic when working in characteristic $2$.

\subsection*{Numerical invariants}
As first immediate consequences we have the following two results.

\begin{Proposition}
  \label{bettiinvariants}
  Let $X$ be a uniruled surface constructed from data $(C,F,\delta)$.
  
  In the Zariski topology, $X$ is homeomorphic  to a
  birationally ruled surface over $C$.
  Moreover,
  $$
  \piet(X) \,\iso\, \piet(C)\,\,\mbox{ and thus }\,\,
  b_1(X)\,=\,b_3(X)\,=\,2g(C).
  $$
  If $g(C)\geq1$, the
  image $B$ of the Albanese map of $X$ is isomorphic 
  to $C$ or $C^{(-1)}$ and the Albanese variety of $X$ 
  is isomorphic to the Jacobian of $B$. 
\end{Proposition}

\prf
By construction, there exists a blow-up $\tilde{S}$ of $C\times\PP^1$
and a finite and purely inseparable morphism $\pi:\tilde{S}\to X$.
The map $\pi$ induces a homeomorphism in the Zariski topology.
By \cite[Th\'eor\`eme IX.4.10]{sga1}, 
the \'etale Betti numbers and
algebraic fundamental groups of $X$ and $\tilde{S}$ coincide.

The Albanese variety of $\tilde{S}$ is the Jacobian of $C$.
Since $\pi:\tilde{S}\to X$ factors over the Frobenius morphism
$F_{\tilde{S}}:\tilde{S}\to\tilde{S}^{(-1)}$, it follows from the
Albanese property that $\alpha(X)$ is a curve that is birational
to $C$ or $C^{(-1)}$.
Hence, the normalisation $B$ of $\alpha(X)$ is isomorphic
to $C$ or $C^{(-1)}$.
By the Albanese property, there exists a map from $\Alb(X)$ to
the Jacobian $\Jac(B)$.
This yields a map from $\alpha(X)$ to $B$, which is an inverse
to the normalisation morphism.
Hence $\alpha(X)$ is isomorphic to $B$.
Applying the Albanese property to $\Alb(X)$ and $\Jac(B)$, it follows
that $\Alb(X)$ is isomorphic to $\Jac(B)$.
\qed

\begin{Proposition}
  \label{canonicaldivisor}
  For data $(C,F,\delta_C+\delta_F)$ the
  quotient $X':=(C\times F)/\delta$ is a normal surface
  and its singularities lie below 
  $(\delta_C)_0\times(\delta_F)_0$ and
  $(\delta_C)_\infty\times(\delta_F)_\infty$.  
  The canonical Weil divisor $K_{X'}$ is $\QQ$-Cartier 
  with self-intersection number
  $$
  K_{X'}^2 \,=\, \frac{2}{p}\,\cdot\, 
  \left( 2g(C)-2 + (p-1)d_C\right)\,\cdot\,
  \left(2g(F)-2 + (p-1)d_F\right)\, ,
  $$
  where $d_C$ denotes the degree of $(\delta_C)_\infty$ and
  $d_F$ denotes the degree of $(\delta_F)_\infty$.
\end{Proposition}

\prf
We already noted that $X'$ is normal and Remark \ref{singularremark}
tells us where to find the singularities of $X'$, which are
isolated points.
A local computation shows that the divisor of $\delta$ is
equal to
\begin{equation}
  \label{divisor}
  (\delta) \,=\, - (\delta_{C})_\infty \cdot F \,-\, 
  (\delta_F)_\infty \cdot C \,.
\end{equation}
Outside the singular locus of $X'$ we have an equality of
Cartier divisors
\begin{equation}
  \label{canonical}
  K_S \,=\, \pi^\ast K_{X'} \,+\, (p-1)\cdot (\delta)\,,
\end{equation}
where $\pi:S\to X'$ is the quotient map.
The Weil divisor $K_{X'}$ corresponds to a reflexive sheaf of 
rank $1$ and so (\ref{canonical}) extends to an equality of Weil 
divisors on the whole of $X'$.
 
Moreover, being reflexive and of rank $1$, the divisor
$K_{X'}$ defines an element in the class group of
every local ring of $X'$.
These class groups are Abelian $p$-torsion groups
since $X'$ is the quotient of a smooth variety by a
$p$-closed derivation, cf. \cite[Chapter IV.17]{fos}.
Hence $pK_{X'}$ is locally principal, i.e. a Cartier 
divisor.
In particular, $K_{X'}$ is a $\QQ$-Cartier divisor.

The assertion on the self-intersection number follows from
(\ref{canonical}) and the projection formula, which we may
use since we are dealing with $\QQ$-Cartier divisors.  
\qed\medskip

The following will be useful later on.
If $\pi:\tilde{S}\to X$ is a finite and purely inseparable
morphism of degree $p$ and height $1$ between smooth varieties
then there exists an exact sequence 
\cite[Corollary 3.4]{ek1}
\begin{equation}
 \label{cotangent}
  0\,\to\,F^\ast\sigma^\ast\Omega_{\tilde{S}/X}\,\to\,
  \pi^\ast\Omega_X^1\,\to\,\Omega_{\tilde{S}}^1
  \,\to\,
  \Omega_{\tilde{S}/X}\,\to\,0\,,
\end{equation}
where $F$ is the $k$-linear
and $\sigma$ is the absolute Frobenius morphism of $\tilde{S}$.

Taking determinants and applying (\ref{canonical}) to $\tilde{S}$
and $X$, we obtain
\begin{equation}
  \label{relative}
\Omega_{\tilde{S}/X}{}^{\otimes (1-p)} \,\iso\, 
\omega_{\tilde{S}}\,\otimes\,\pi^\ast\omega_X{}^\vee \,\iso\,
\OO_{\tilde{S}}( (p-1)\cdot(\widetilde{\delta}))\,,
\end{equation}
where $\tilde{\delta}$ is the rational vector field defining $\pi$.

\subsection*{Singular fibrations}
We will now show that
our surfaces are endowed with two fibrations, both of which are
usually not generically smooth.
Together with the results of Section \ref{curvesection}, 
we use the following proposition
to construct fibrations with prescribed singularities 
of the geometric generic fibre.

Since regularity, and hence also normality,
are not stable under inseparable field extensions, 
we define the {\em arithmetic genus}
$p_a$ of an irreducible curve $C$, which is defined over a possibly 
non-perfect ground field $L$, to be $1-\chi(\OO_C)$.
For such a curve, we define its {\em geometric genus} $g$
to be the genus of the normalisation of 
$C\otimes_L \bar{L}$, where $\bar{L}$ is an algebraic
closure of $L$.

\begin{Proposition}
  \label{cuspfib}
  Let $X$ be a surface constructed from data 
  $(C,F,\delta=\delta_C+\delta_F)$ and denote by $d_F$
  the degree of the divisor of poles $(\delta_F)_\infty$
  of $\delta_F$.
  We assume that $\delta_C$ is non-trivial.
  
  Then there exists a fibration
  $f\,:\, X \,\to\, C^{(-1)}$.
  Its generic fibre $F_\eta'$ is a regular curve over $k(C)^p$
  of arithmetic genus
  $$
  p_a(F_\eta')\,=\, g(F) \,+\, {\textstyle\frac{p-1}{2}}\,\cdot\, d_F\,.
  $$
  The normalisation of the geometric generic fibre 
  is isomorphic to $F\otimes_k \overline{k(C)^p}$ 
  and its singular points 
  are cusps
  lying below the points where $\delta_F$ has a pole.
  More precisely, if $\delta_F$ has a pole of order 
  $m$ in a point $P\in F$,
  then the arithmetic genus of the singularity of the geometric
  generic fibre at the cusp lying below $P$ is equal to
  $\frac{p-1}{2}\cdot m$.
\end{Proposition}

\prf
By definition, there exists a blow-up $\tilde{S}$ of $S:=C\times F$
and a finite inseparable morphism $\pi:\tilde{S}\to X$ of height $1$.
Factoring the Frobenius morphism, we obtain a morphism 
$\varpi:X\to \tilde{S}^{(-1)}$.
The latter surface has a morphism onto $S^{(-1)}$ and a projection
onto its factor $C^{(-1)}$.
Composing, we obtain a morphism $f$ from $X$ onto $C^{(-1)}$.

The generic fibre of the fibration of $\tilde{S}^{-1}$ over $C^{(-1)}$ 
is $F^{(-1)}$.
Since $\varpi$ is purely inseparable, the generic fibre of $f$ is
homeomorphic to $F^{(-1)}$.
As $\varpi$ has height $1$, the normalisation of the
geometric generic fibre of $f$ is isomorphic to 
$F\otimes_k \overline{k(C)^p}$.
The normalisation map is easily seen to be a homeomorphism, 
and so the singularities of the geometric generic fibre are 
unibranch, i.e. cusps. 

The local rings of the generic fibre are 
localisations of the local rings of the surface $X$, 
which is normal.
Thus, these rings are normal and
being $1$-dimensional, they are regular.
Hence the generic fibre is regular.

To compute the arithmetic genus of the geometric generic fibre
$A:=F_\eta'\otimes_k \overline{k(C)^p}$,
we use the adjunction formula
$2p_a(A)-2=K_X\cdot F_\eta'+ F_\eta'^2$.
Since $F_\eta'$ is a fibre, its self-intersection number 
is zero.
We assumed $\delta_C$ to be non-trivial.
Hence $F$ is not an integral curve for $\delta$ and
so $\pi_\ast F=F_\eta'$
by \cite[\S1 Proposition 1]{rs}. 
Then the projection formula yields
$F_\eta'\cdot K_X=F\cdot \pi^\ast K_X$.
Since we are dealing with the generic fibre we can ignore
contributions coming from the exceptional divisors on
$\tilde{S}$.
Thus we compute the intersection numbers on the singular
surface $X'$ and obtain the asserted formula for $p_a(A)$
using formulae (\ref{divisor}) and (\ref{canonical}).

Performing the previous calculations 
locally, we see that the singularities of the 
(geometric) generic fibre correspond to those 
points where $\delta_F$ has a pole.
\qed

\begin{Remark}
This result illustrates Tate's theorem \cite{t} 
that the genus change of a curve in
inseparable field extensions is divisible by $\frac{p-1}{2}$.
\end{Remark}

In the analysis of hyperelliptic surfaces,
one of the two (quasi-) elliptic fibrations comes from 
the Albanese map.
Joining the proofs of Proposition \ref{bettiinvariants} 
and Proposition \ref{cuspfib} together we obtain

\begin{Corollary}  
  If $F$ is rational and $g(C)\geq1$, then this fibration
  coincides with the Stein factorisation of the Albanese morphism
  of $X$.\qed
\end{Corollary}

\section{Hodge theory in characteristic $2$}

To get finer invariants, we work in characteristic $2$,
where some of the computations become easier.
Given $\pi:\tilde{S}\to X$, there exists a unique
morphism $\varpi:X\to\tilde{S}^{(-1)}$ such that
$\varpi\circ\pi$ is equal to the Frobenius map
of $\tilde{S}$.
Since $\varpi$ is flat of degree $2$, there exists
a line bundle $\cal L$ and a short exact sequence
\begin{equation}
  \label{deflinebundle}
  0\,\to\,\OO_{\tilde{S}^{(-1)}}\,\to\,\varpi_\ast\OO_X\,\to\,
  {\cal L}^\vee\,\to\,0\,.
\end{equation}
By \cite[Proposition I.1.11]{ek2}, the map $\varpi$ has the
structure of a torsor under the finite flat
group scheme $\alpha_{\cal L}$.
We recall that $\alpha_{\cal L}$ is defined to be the
kernel of the Frobenius map from $\cal L$ to ${\cal L}^p$
where we regard these line bundles as group schemes in the
flat topology on $\tilde{S}^{(-1)}$.

A nice feature is that a possibly singular
$\alpha_{\cal L}$-torsor $X$ over a smooth base 
$\tilde{S}^{(-1)}$ is automatically Gorenstein.
More precisely
\begin{equation}
  \label{dualisingsheaf}
  \omega_X \,\iso\, \varpi^\ast
  (\omega_{\tilde{S}^{(-1)}}\otimes{\cal L}^{\otimes(p-1)})\,, 
\end{equation}
where $\omega_X$ denotes the dualising sheaf of $X$, cf.
\cite[Proposition I.1.7]{ek2}.

\subsection*{Hodge numbers}
With these preliminary remarks we determine Hodge numbers.

\begin{Proposition}
 \label{hodgeinvariants}
  Let $X$ be a uniruled surface constructed from data 
  $(C,F,\delta=\delta_C+\delta_F)$
  in characteristic $2$.
  We let $d_C$ and $d_F$ be the degrees of the divisors of
  poles $(\delta_C)_\infty$ and $(\delta_F)_\infty$, 
  respectively.
  
  The numerical equivalence class $L$ of 
  the line bundle $\cal L$ is given by
  \begin{equation}
    \label{linebundle}
    -L \,\num\, \left(-\frac{d_F}{2}-1\right)\, C \,+\, 
    \left(-\frac{d_C}{2}+g(C)-1\right)\, F \,+\,\sum_i a_i\, E_i\,,
  \end{equation}
  for some integers $a_i>0$,
  where the $E_i$'s are the exceptional divisors of the blow-up
  $\tilde{S}\to S$.
  Then,
  $$\begin{array}{lcl}
    \blowup \chi(\OO_X) &=& 2(1-g(C)) \,+\, {\textstyle \frac{1}{2}} (L^2+K_{\tilde{S}}L)  \\
    \blowup h^{01}(X) \,=\, h^1(\OO_X) &=& g(C) \,+\, h^1( {\cal L}^\vee ) \\
    \blowup h^{02}(X) \,=\, p_g(X) &=& h^2({\cal L}^\vee)\,.
    \end{array}$$
  In particular, 
  the Picard scheme of $X$ is reduced if and only if $h^1({\cal L}^\vee)=0$.
\end{Proposition}

\prf
The formula for $\chi(\OO_X)$ follows from (\ref{deflinebundle}) 
and Riemann-Roch on $\tilde{S}^{(-1)}$.
The long exact sequence in cohomology applied to 
(\ref{deflinebundle}) together with the fact that 
$h^2(\OO_{\tilde{S}^{(-1)}})$ vanishes (we assumed 
$F$ to be rational), yields the formulae for 
$h^{01}$ and $h^{02}$.
The Albanese variety is $g(C)$-dimensional by 
Proposition \ref{bettiinvariants}.
Thus, $\Pic(X)$ is reduced iff
$h^1(\OO_X)=g(C)$, i.e. iff $h^1({\cal L}^\vee)=0$.

To determine $\cal L$, we run through the arguments given
in \cite[Section I.1]{ek2}.
The morphism $\varpi$ is an $\alpha_{{\cal L}}$-torsor.
This defines an embedding of ${\cal L}^{\otimes (-p)}$
into $\Omega_{\tilde{S}^{(-1)}}^1$.
The annihilator 
${\cal M}\subseteq\Theta_{\tilde{S}^{(-1)}}^1$
is the $1$-foliation defining the morphism
$\pi^{(-1)}:\tilde{S}^{(-1)}\to X^{(-1)}$ and
we obtain a short exact sequence
\begin{equation}
\label{foliationseq}
0\,\to\,{\cal M}\,\to\,\Theta_{\tilde{S}^{(-1)}}^1\,\to\,
{\cal L}^{\otimes p}\,\to\,0\,.
\end{equation}
The line bundle $\cal M$ is isomorphic to
$\Omega_{\tilde{S}^{(-1)}/X^{(-1)}}^\vee$.
Taking determinants we obtain for $p=2$ and together
with (\ref{cotangent}) and (\ref{relative}) an equality
\begin{equation}
\label{determinant}
  {\cal L}^{\otimes (-2)} \,\iso\,
  \omega_{\tilde{S}^{(-1)}} \,\otimes\, 
  \OO_{\tilde{S}^{(-1)}}((\tilde{\delta}))\,.
\end{equation}
The divisor of $\delta$ is given by (\ref{divisor}).
Resolving the singularities of $\delta$ via a blow-up
$\tilde{S}\to S$, the divisor $\tilde{\delta}$ of
the pull-back of $\delta$ is given by
$$
(\tilde{\delta}) \,=\, -(\delta_C)_\infty\cdot F\,-\, 
(\delta_F)_\infty\cdot C \,+\, \sum_i b_i E_i
$$
for some non-negative integers $b_i$ as explained in the
proof of \cite[Proposition 2.6]{hir}.
From this and (\ref{determinant}), formula (\ref{linebundle}) follows 
immediately.
\qed

\begin{Proposition}
  \label{hodgeinvariants2}
  Under the assumptions and notations of 
  Proposition \ref{hodgeinvariants} we have 
  $$\begin{array}{lclcl}
    h^{10}(X) \,=\, h^0(\Omega_X^1) &=& 
    h^0(\tilde{S},\,\varpi_\ast\varpi^\ast\OO_{\tilde{S}^{(-1)}}(-\tilde{\delta}))
    &\geq& g(C)\,.
    \end{array}$$ 
  There exists a short exact sequence
  $$
  0\,\to\,H^0(\tilde{S},\OO_{\tilde{S}}(-\tilde{\delta}))\,\to\,
  H^0(X,\Omega_X^1)\,\stackrel{d_X}{\to}\,H^0(X,\,\Omega_X^2)\,,
  $$
  where $d_X$ denotes the differential.
\end{Proposition}

\prf
We rewrite the exact sequence (\ref{cotangent}) as
\begin{equation}
  \label{frobcotangent}
  0\,\to\,F^\ast\sigma^\ast\Omega_{X/\tilde{S}^{(-1)}}
  \,\stackrel{\beta}{\to}\,
  \varpi^\ast\Omega_{\tilde{S}^{(-1)}}^1\,\to\,\Omega_X^1
  \,\to\,
  \Omega_{X/\tilde{S}^{(-1)}}\,\to\,0\,.
\end{equation}
From (\ref{relative}) and (\ref{dualisingsheaf}) we obtain
$\Omega_{X/\tilde{S}^{(-1)}}\iso\varpi^{\ast}({\cal L}^\vee)$.
Using (\ref{deflinebundle}), we see that
$\varpi_\ast\Omega_{X/\tilde{S}^{(-1)}}$ is an 
extension of ${\cal L}^{\otimes (-2)}$ by ${\cal L}^\vee$.
It is not difficult to see that neither of these line bundles 
has global sections.
Hence $\Omega_{X/\tilde{S}^{(-1)}}$ has no non-trivial global
sections.

We take determinants in (\ref{frobcotangent})
and use (\ref{foliationseq}) as well as
(\ref{determinant}) to
conclude that the cokernel of $\beta$ is the pullback
$\varpi^\ast{\cal E}$ of the cokernel $\cal E$ of
\begin{equation}
\label{ses}
0\,\to\,{\cal L}^{\otimes (-2)}\,\to\,\Omega_{\tilde{S}^{(-1)}}^1
\,\to\,\OO_{\tilde{S}^{(-1)}}(-\tilde{\delta})\,\to\,0\,.
\end{equation}
By what we have already shown, we know that
$h^0(X,\Omega_X^1)=h^0(X,\varpi^\ast{\cal E})=
h^0(\tilde{S},\varpi_\ast\varpi^\ast{\cal E})$.
By the projection formula and (\ref{deflinebundle}),
$\varpi_\ast\varpi^\ast{\cal E}$ is an extension
$$
0\,\to\,\OO_{\tilde{S}^{(-1)}}(-\tilde{\delta})
\,\to\,\varpi_\ast\varpi^\ast{\cal E}\,\to\,
\OO_{\tilde{S}^{(-1)}}(-\tilde{\delta})\otimes{\cal L}^{\vee}
\,\to\,0\,.
$$
By (\ref{dualisingsheaf}), the pull-back $\varpi^\ast$ of
the cokernel is isomorphic to $\omega_X$.
That the induced morphism on global sections from
$H^0(X,\Omega_X^1)\iso H^0(\tilde{S}^{(-1)},\varpi_\ast\varpi^\ast{\cal E})$
to
$H^0(\tilde{S}^{(-1)},{\cal E}\otimes{\cal L}^\vee)\iso H^0(X,\omega_X)$
coincides with the differential map follows from
\cite[Theorem 3.1]{hir} applied to $\tilde{S}$ and $\tilde{\delta}$.
\qed

\subsection*{Rational singularities}

The cohomology of line bundles on $C\times F$ is easily computed.
However, in order to obtain the Hodge invariants of a surface
$X$ constructed from data $(C,F,\delta)$, we have to desingularise
$X'\to (C\times F)^{(-1)}$ and in order to use
Proposition \ref{hodgeinvariants} and 
Proposition \ref{hodgeinvariants2}
we have to compute the cohomology of line bundles on a blow-up
of $(C\times F)^{(-1)}$, which is usually more complicated.
Things become easier if we assume that $X'$ has at worst rational
singularities.

A non-trivial vector field $\delta_C$ on a curve of genus 
$g(C)$ has $d_C:=\deg(\delta_C)_\infty\geq 2-2g(C)$.
The Hodge invariants become easier to compute if we assume
that equality does not hold:

\begin{Proposition}
  \label{hodgeinvariants3}
  Let $X$ be a uniruled surface constructed from data
  $(C,F,\delta)$.
  Assume that the singular quotient $X':=S/\delta$ has 
  at worst rational singularities
  and $d_C>2g(C)-2$.
  Then
  \begin{enumerate}
  \item $h^1(X,\,\OO_X)\,=\,g(C)$ and 
  \item $h^0(X,\,\Omega_X^1)\,=\,g(C)$.
  \end{enumerate}
  In particular,
  the Picard scheme of $X$ is reduced,
  all global $1$-forms on $X$ are 
  pull-backs of global $1$-forms from $\Alb(X)$
  and all global $1$-forms on $X$
  are $d$-closed.
\end{Proposition}

\prf
We consider the finite flat morphism
$\varrho:X'\to S^{(-1)}$.
We define a line bundle $\cal N$ on $S^{(-1)}$
associated to $\varrho$ as in (\ref{deflinebundle}).
Then ${\cal N}^\vee$ is 
numerically equivalent to
$$
{\cal N}^{\vee}\,\num\,
\OO_{S^{(-1)}}\left( 
(-\frac{1}{2}d_F-1)C \,+\, (-\frac{1}{2}d_C+g(C)-1)F
\right)\,.
$$
We assumed $d_C>2g(C)-2$ and so the K\"unneth formula
yields $h^1({\cal N}^{\otimes(-i)})=0$ for $i\geq1$.

Since $X'$ has at worst rational singularities, we
have an equality 
$h^1(\OO_X)=h^1(\OO_{X'})=
g(C)+h^1(S^{(-1)},{\cal N}^\vee)$, which is equal to
$g(C)$ by the previous paragraph.
Hence the Picard scheme of $X$ is reduced.

Since $X'$ is Gorenstein and we assumed that it has rational
singularities it assumption it has Du~Val singularities 
only.
On the other hand, the resolution of singularities $p:X\to X'$
is minimal by \cite[Lemma 2.1]{sbfol} and we
obtain $p^\ast\omega_{X'}\iso\omega_X$.
For the dualising sheaves of $X$ and $X'$ we have
the formulae
$\omega_X\iso\varpi^\ast(\omega_{\tilde{S}}\otimes{\cal L})$
and
$\omega_{X'}\iso\varrho^\ast(\omega_{S}\otimes{\cal N})$.
From this we deduce an isomorphism
\begin{eqnarray}
\label{related}
{\cal L} &\iso& f^\ast({\cal N}\otimes\omega_S)\,\otimes\,\omega_{\tilde{S}}^\vee\,,
\end{eqnarray}
where $f:\tilde{S}\to S$ denotes the blow-up needed to resolve the
isolated singularities of $\delta$.

By Proposition \ref{hodgeinvariants2}, the global sections of 
$\OO_{\tilde{S}^{(-1)}}{(-\tilde{\delta})}$ correspond to the global
$d$-closed $1$-forms on $X$.
Pushing forward (\ref{ses}) to $S^{(-1)}$, 
we obtain a long exact sequence
$$
\begin{array}{ll}
0\,\to & {\cal N}^{(-2)}\,\to\,\Omega_{S^{(-1)}}^1\,\to\,
f_\ast\OO_{\tilde{S}^{(-1)}}{(-\tilde{\delta})}\,\to\,\\
&R^1f_\ast{\cal N}^{(-2)}\,\to\,R^1f_\ast\Omega_{S^{(-1)}}^1\,\to\,
R^1f_\ast\OO_{\tilde{S}^{(-1)}}{(-\tilde{\delta})}\,\to\,0\,.
\end{array}
$$
It is not difficult to see that 
$R^1f_\ast\OO_{\tilde{S}^{(-1)}}{(-\tilde{\delta})}=0$.
If $f:\tilde{S}\to S$ is a sequence of $r$ blow-ups along closed points
then by an induction on the number of blow-ups and
using (\ref{related}) we conclude that
$R^1f_\ast{\cal N}^{(-2)}$ and $R^1f_\ast\Omega_{S^{(-1)}}^1$
both are Artin algebras of length $r$.
In particular, we obtain a short exact sequence
\begin{equation}
\label{dclosedforms}
0\,\to\,{\cal N}^{(-2)}\,\to\,\Omega_{S^{(-1)}}^1\,\to\,
f_\ast\OO_{\tilde{S}^{(-1)}}{(-\tilde{\delta})}\,\to\,0\,.
\end{equation}
Taking cohomology and noting that $H^0({\cal N}^{(-2)})=H^1({\cal N}^{(-2)})=0$ 
we see that the space of $d$-closed global $1$-forms 
on $X$ is $g(C)$-dimensional.

Let $x,y$ be local coordinates on $\tilde{S}^{(-1)}$ such that $y=0$ defines
a component $E'$ of the exceptional divisor of the blow-up $f$.
We also assume that $E'$ is the only possible component of the divisor of
$\tilde{\delta}$ through the point $(x,y)=(0,0)$. 
We denote by $e$ the order of pole
$\OO_{\tilde{S}^{(-1)}}(-\tilde{\delta})$ has along $E'$,
a number which can be read off from (\ref{related}).
If $E'$ is an integral curve for the foliation
$\tilde{\delta}$ then 
$\varpi^\ast\OO_{\tilde{S}^{(-1)}}(-\tilde{\delta})$
is locally generated by
$(\varpi^\ast y)^e \, d (\varpi^\ast y)$
as a subsheaf of $\Omega_X^1$.
Then, $\varpi^\ast y$ and a square root $\sqrt{\varpi^\ast x}$
form a system of local coordinates on $X$.
If $E'$ is not an integral curve for $\tilde{\delta}$ then
$\varpi^\ast\OO_{\tilde{S}^{(-1)}}(-\tilde{\delta})$
is locally generated by
$(\varpi^\ast y)^e \, d (\varpi^\ast x)$.
In this case, $\varpi^\ast x$ and a square root
$\sqrt{\varpi^\ast y}$ form a system of local coordinates on $X$.

We have already seen in the proof of Proposition \ref{hodgeinvariants2}
that all global $1$-forms are global sections of
$\varpi^\ast\OO_{\tilde{S}^{(-1)}}(-\tilde{\delta})$.
Whence the derivative of a global $1$-form is can be computed
on forms of the form a regular function times 
$(\varpi^\ast y)^e \, d (\varpi^\ast x)$ and
$(\varpi^\ast y)^e \, d (\varpi^\ast y)$, respectively.
Using the local coordinates on $X$ above, we
see that differentiating such a $1$-form we get a $2$-form, which
has a zero of order at least $e$ along $\varpi^\ast(y)=0$.
On the other hand, $X'$ has Du~Val singularities only, 
which have no adjunction condition, i.e. $\Omega_X^2$ is 
locally around $(x,y)=(0,0)$ generated by 
$d(\varpi^\ast y)\wedge d\sqrt{\varpi^\ast x}$
if $E'$ is an integral curve and
$d(\sqrt{\varpi^\ast y})\wedge d(\varpi^\ast x)$
if $E'$ is not an integral curve for $\tilde{\delta}$.

We define 
${\cal F}\,:=\,
\omega_{S^{(-1)}}\otimes\omega^\vee_{\tilde{S}^{(-1)}}$, which we
consider as subsheaf of $\OO_{\tilde{S}^{(-1)}}$.
By (\ref{related}) and our local computations, we see that
$d\Omega_X^1$ is a subsheaf of
$\varpi^{-1}({\cal F})\cdot\Omega_X^2\,=\,\varpi^\ast{\cal F}\otimes\Omega_X^2$.
Taking global sections we get
$$
H^0(\Omega_X^1)\,\stackrel{d_X}{\to}\,
H^0( \varpi^\ast({\cal F})\,\otimes\,\Omega_X^2) 
\,\subseteq\,H^0(\Omega_X^2)\,.
$$

In order to show that all global $1$-forms are $d$-closed it is enough
to show that the dimension of the space in the middle is zero.
Pushing the sheaf forward to $S^{(-1)}$ and using
(\ref{dualisingsheaf}) we obtain an extension
$$
0\,\to\,
{\cal F}\otimes
(\omega_{\tilde{S}^{(-1)}}\otimes{\cal L})
\,\to\,
\varpi_\ast
(\varpi^\ast({\cal F})\,\otimes\,\Omega_X^2)
\,\to\,
{\cal F}\otimes\omega_{\tilde{S}^{(-1)}}
\,\to\,0\,.
$$
The sheaf on the right has only trivial global sections and so
we have to compute those on the left.
With (\ref{dclosedforms})
we obtain a short exact sequence
\begin{equation}
\label{dclosedforms2}
0\,\to\,{\cal N}^{(-3)}\,\to\,\Omega_{S^{(-1)}}^1\otimes{\cal N}^\vee\,\to\,
f_\ast\left( {\cal F}
\,\otimes\,(\omega_{\tilde{S}^{(-1)}}\otimes{\cal L})
\right)
\,\to\,0\,.
\end{equation}
Using the vanishing of $H^1({\cal N}^{(-3)})$ and 
$H^0(\Omega_S^1\otimes{\cal N}^\vee)$, it finally follows that
all global $1$-forms on $X$ are $d$-closed.

The Albanese variety of $X$ is 
$g(C)$-dimensional by Proposition \ref{bettiinvariants}.
Igusa's theorem \cite{ig} states that the pull-back of a non-trivial
global $1$-form on $\Alb(X)$ to $X$ remains non-trivial.
This implies $h^0(\Omega_X^1)\geq g(C)$.
If equality holds then every global $1$-from on 
$X$ is the pull-back of a global $1$-form on $\Alb(X)$.
\qed

\begin{Remark}
  Being Gorenstein and having rational singularities, $X'$ has at
  worst Du~Val singularities.
  It is not difficult to see that these can be of type
  $A_1$, $D_{2n}$, $E_7$ and $E_8$ only, cf.
  \cite[Lemma 2.1]{sbfol}.
\end{Remark}

\subsection*{The Fr\"olicher spectral sequence}
For a smooth variety $X$, Hodge and de~Rham cohomology
are related by the so-called
{\it Fr\"olicher spectral sequence}
$$
E^{ij}_1\,:=\,H^j(X,\,\Omega^i_{X/k}) \,\Rightarrow\, \HdR{i+j}(X/k) .
$$
It follows from classical Hodge theory that this
spectral sequence degenerates at $E_1$-level for
K\"ahler manifolds.
For curves and complex surfaces we even have degeneration 
at $E_1$-level without the K\"ahler assumption, cf.
\cite[Chapter IV.2]{bhpv}.
However, for surfaces over arbitrary fields, there is no reason for 
this spectral sequence to degenerate at $E_1$-level.

\begin{Theorem}
  \label{hodgederham}
  Let $X$ be a uniruled surface
  constructed from data $(C,F,\delta)$ in characteristic $2$
  such that $(C\times F)/\delta$ has at worst rational 
  singularities.
  Then the crystalline cohomology of $X$ is torsion-free and its
  Fr\"olicher spectral sequence degenerates at
  $E_1$-level.
\end{Theorem}

\prf
Since $\Hcris{1}(X/W)$ is a torsion-free $W$-module,
we know from Proposition \ref{bettiinvariants} that
its rank is equal to $2g:=2g(C)$.

We consider the universal coefficient formula
\begin{equation}
\label{universalcoefficient}
0\,\to\,\Hcris{i}(X/W)\otimes_{W}k\,\to\,\HdR{i}(X/k)\,\to\,
{\rm Tor}_1^W(\Hcris{i+1}(X/W),\,k)\,\to\,0\,.
\end{equation}
Thus, $\HdR{1}(X/k)$ is at least $2g$-dimensional.
The existence of the Fr\"olicher spectral sequence
already implies the inequality 
$\hdR{1}\leq h^{01}+h^{10}$.
Since $h^{01}=h^{10}$ by Proposition \ref{hodgeinvariants}
and Proposition \ref{hodgeinvariants3} we have equality.
Because $\Hcris{4}(X/W)$ is torsion-free, the universal 
coefficient formula (\ref{universalcoefficient}) for $i=3$
yields $\hdR{3}=2g$.
By Serre duality, we obtain 
$\hdR{3}=2g=h^{01}+h^{10}=h^{12}+h^{21}$.
Again, it follows already from the existence of the
Fr\"olicher spectral sequence that the sum over
the $(-1)^i\hdR{i}$ is equal to the sum 
$(-1)^{i+j}h^{ij}$.
Since we are working with surfaces,
this implies $\hdR{2}=h^{02}+h^{11}+h^{20}$.
Hence $\hdR{n}$ is equal to the sum over all 
$h^{ij}$ with $i+j=n$ for all $n$.
This implies that the Fr\"olicher spectral sequence degenerates at
$E_1$-level.

Plugging this into (\ref{universalcoefficient}), we see that
${\rm Tor}_1^W(\Hcris{i}(X/W),\,k)=0$ for all $i$ which implies
that the crystalline cohomology of $X$ is torsion-free. 
\qed

\subsection*{The slope spectral sequence}
In \cite[Section II.3]{ill}, Illusie constructs a
a spectral sequence from Hodge-Witt cohomology 
to crystalline cohomology
$$
E^{ij}_1\,:=\,H^j(X,\,W\Omega^i_X) \,\Rightarrow\, \Hcris{i+j}(X/W) .
$$
Modulo torsion, this sequence always degenerates 
at $E_1$-level, cf. \cite[Th\'eor\`eme II.3.2]{ill}.
In general, degeneracy at $E_1$-level
is equivalent to the torsion subgroups of the $H^j(W\Omega_X^i)$'s 
being finitely generated $W$-modules.
For surfaces, Nygaard has shown that this is
equivalent to the finite generation of $H^2(W\OO_X)$, cf.
\cite[Corollaire II.3.14]{ill}.

If the slope spectral sequence degenerates at $E_1$-level
then the crystalline cohomology decomposes
as a direct sum of the Hodge-Witt cohomology groups
and the variety is said to be of
{\it Hodge-Witt type}.

\begin{Theorem}
  \label{slope}
  Let $X$ be a uniruled surface constructed from data $(C,F,\delta)$
  in characteristic $2$ such that $(C\times F)/\delta$ has at worst
  rational singularities.
  If $\chi(\OO_X)>1-g(C)$ then the slope spectral sequence does not
  degenerate at $E_1$-level.
\end{Theorem}

\prf
By the previous result we know that the crystalline cohomology of
$X$ is torsion-free.
Since $\tilde{S}$ is birationally ruled, 
$\Hcris{2}(X/W)\otimes_W K$ is pure of slope $1$ and so
$\Hcris{2}(X/W)\otimes_W K$ is pure of slope $1$.
Being torsion-free, already $\Hcris{2}(X/W)$ 
is pure of slope $1$.

Suppose that the slope spectral sequence degenerates.
Then we have a Hodge-Witt decomposition of the crystalline 
cohomology of $X$, cf. \cite[Th\'eor\`eme IV.4.5]{illray}.
In particular, $H^2(W\OO_X)$ can be identified with the part of
$\Hcris{2}(X/W)$ that has slope strictly less than $1$.
Since $\Hcris{2}(X/W)$ is pure of slope $1$, we see that
$H^2(W\OO_X)$ is zero.

For every $n\geq1$, the Verschiebung $V$ induces a short exact sequence
$$
0\,\to\, VW_{n-1}\OO_X \,\to\, W_n\OO_X \,\to\, \OO_X\,\to\,0\,.
$$
Taking cohomology, passing to the inverse limit and noting that 
$H^3(VW_{n-1}\OO_X)$ vanishes for all $n$, we obtain a surjective
homomorphism of $W$-modules from $H^2(W\OO_X)$ onto $H^2(\OO_X)$.

However, $H^2(\OO_X)\neq0$ since we assumed $\chi(\OO_X)>1-g$
whereas $H^2(W\OO_X)$ is zero.
This contradiction shows that the slope spectral sequence does
not degenerate at $E_1$-level.\qed\medskip

For a smooth variety $X$ the sheaf $B\Omega_X^i$ is defined to
be the image of $d:\Omega_X^{i-1}\to\Omega_X^i$.
Then $X$ is called {\em ordinary} if
$H^j(X,B\Omega_X^i)=0$ for all $i,j$.
By \cite[Th\'eor\`eme IV.4.13]{illray}, the slope
spectral sequence of an ordinary variety degenerates
at $E_1$ and thus we obtain 

\begin{Corollary}
  A uniruled surface constructed from data $(C,F,\delta)$
  as in Theorem \ref{slope}
  is not ordinary.
  \qed
\end{Corollary}

\begin{Remark}
 If a surface $X$ is fibred over some curve $C$ such that
 the generic fibre is a {\em smooth} rational curve, then $X$
 is birationally ruled and such a surface is ordinary if and only
 if the curve $C$ is ordinary.
\end{Remark}

\section{Arithmetic observations}

\subsection*{Artin invariants}
A basic invariant of a surface is its Picard number $\rho$, 
i.e. the rank of its N\'eron--Severi group.
By the Igusa--Severi inequality, we always have $\rho\leq b_2$
and a surface is called {\it supersingular in the sense of Shioda}
if equality holds.
When Artin studied supersingular K3 surfaces \cite{ar}, he observed
that the discriminant of the intersection form on the N\'eron--Severi
group is always negative and an even power of the characteristic 
of the ground field.

\begin{Proposition}
  Let $X$ be a uniruled surface constructed from data $(C,F,\delta)$.
  
  Then $X$ is supersingular in the sense of Shioda.
  The discriminant of its N\'eron--Severi lattice is
  $$
  {\rm disc}\,\NS(X)\,=\, -\,p^{2\sigma}\,,
  $$
  where $\sigma$ is a non-negative integer, i.e. there exists
  an Artin invariant $\sigma$ for such surfaces.
  
  It is bounded by $\sigma\leq b_2/2$.
  If the crystalline cohomology of $X$ is torsion-free
  then $\sigma\geq p_g$.
\end{Proposition}

\prf
There exists a finite purely inseparable map 
$\pi$ from a birationally ruled surface $\tilde{S}$ onto $X$, 
where $S$ is ruled over the curve $C$.
The intersection form on $NS(\tilde{S})$ is unimodular.
By the Hodge index theorem, its signature is 
equal to $(1,\rho(\tilde{S})-1)$.
Hence its discriminant is equal to
$(-1)^{\rho(\tilde{S})-1}$.

Since $\pi$ is inseparable, $\pi^\ast NS(X)$ is a subgroup of finite 
index in $NS(\tilde{S})$.
In particular, we have $\rho(X)=\rho(\tilde{S})=:\rho$, i.e. $X$ is
supersingular in the sense of Shioda.
As $\pi$ has degree $p$, this index is equal 
to $p^n$ for some integer $n$.
Hence the discriminant of $\pi^\ast NS(X)$ is equal to 
$(-1)^{\rho-1}p^{2n}$.
Using the projection formula 
$\pi^\ast A\cdot\pi^\ast B = p A\cdot B$, it follows that the
discriminant of $NS(X)$ is equal to 
$(-1)^{\rho-1}p^{2n-\rho}$.
So once we have shown that $2n-\rho$ is even it follows that
$\rho$ is even and that the sign of the discriminant is $-1$.

The first Chern class provides us with an injective homomorphism
$$
c_1\,:\, NS(X)\otimes W \,\into\, \Hcris{2}(X/W)
$$
where the intersection pairing on the left hand side coincides
with the cup-product on the right hand side,
cf. \cite[Remarque II.5.21.4]{ill}.
By Poincar\'e duality, the cup-product on $\Hcris{2}(X/W)$ 
(modulo torsion) is unimodular. 
Both sides are free of rank $\rho=b_2$ and so the image of
$c_1(NS(X))$ in $\Hcris{2}(X/W)$ has finite index.
This index is a $p$-power and it follows that the 
discriminant of $NS(X)\otimes W$ is an even $p$-power.

Let $F:\tilde{S}\to\tilde{S}^{(-1)}$ be the Frobenius morphism.
Then $F^\ast(NS(\tilde{S}^{(-1)}))$ is a sublattice
of $\pi^\ast NS(X)$ since $F$ factors over $\pi$.
The index of $F^\ast (NS(\tilde{S}^{(-1)}))$ in $NS(\tilde{S})$
is equal to $p^{\rho}$. 
Thus, the index $p^{n}$
of $\pi^\ast NS(X)$ in $NS(\tilde{S})$ divides $p^{\rho}$.
This implies $n\leq\rho$, and so $2\sigma=2n-\rho\leq\rho=b_2$.

If the crystalline cohomology of $X$ is torsion-free then 
$\Pic(X)$ is reduced by \cite[Proposition II.5.16]{ill}.
Then we can argue as in \cite[Remarque II.5.21]{ill} to
conclude $\sigma\geq p_g(X)$.
\qed\medskip

It is known that the N\'eron--Severi lattice of a $K3$ surface
is even and unimodular.
More precisely, it is a sublattice of $3H\oplus(-2)E_8$,
where $H$ is a hyperbolic plane. 
By a result of Rudakov and \v{S}afarevi\v{c}, the Artin 
invariant  determines the intersection form of a supersingular K3 surface
up to isomorphism.
In our case, the situation is more complicated.
In some examples, we have to resolve an elliptic $(19)_0$-singularity 
to obtain our surface. Then there is a curve with 
self-intersection $-3$ and so the intersection form is not even.

\subsection*{The Artin--Tate conjecture}
Let $X$ be a surface over the finite field $k$ with $q$ elements
and let $Z(X,t)$ be its zeta function.
By Deligne's proof of the Weil conjectures it is known that
there exist polynomials $P_1$, $P_2$ and $P_3$ of degrees
equal to the Betti numbers $b_1$, $b_2$ and $b_3$ of $X$ such that
$$
 Z(X,\,t)\,=\, \frac{P_1(X,\,t)\,\cdot\, P_3(X,\,t)}%
 {(1-q\, t)\,\cdot\,P_2(X,\,t)\,\cdot\,(1-q^2\, t)}\,.
$$
We denote by $\rho(X)$ the Picard number and by $\Br(X)$ 
the Brauer group of $X$.
We define $\alpha(X):=\chi(\OO_X)-1+b_1(X)/2$.
The Artin--Tate conjecture states that
$$
 P_2(X,\,q^{-s}) \,\sim\, 
 (-1)^{\rho(X)-1}\,\cdot\,
 \frac{ {\rm disc}\,\NS(X)\,\cdot\,|\Br(X)|}
 {|NS(X)_{\rm tors}|\,\cdot\,q^{\alpha(X)}}
 \,\cdot\,
 (1-q^{1-s})^{\rho(X)}\mbox{ as } s\to1\,.
$$
\begin{Proposition}
  Let $X$ be a surface constructed from data $(C,F,\delta)$,
  where $C$, $F$ and $\delta$ are defined over the finite field $k$
  with $q$ elements.
  
  Then $X$ is defined over $k$, the Artin--Tate conjecture holds
  for $X$, and we have an equality
  $$
   {\rm disc}\,\NS(X) \,\cdot\, |\Br(X)| \,=\, q^{\alpha(X)\,-\,g(C)\cdot g(F)}\,
   \cdot\, |NS(X)_{\rm tors}|\,.
  $$    
\end{Proposition}

\prf
  To construct $X$ we had to find a blow-up $\tilde{S}$ of $S=C\times F$ 
  and quotiented by a vector field to obtain $\pi:\tilde{S}\to X$.
  All this can be done over $k$ and so $X$ is defined over $k$.
  
  Since the Tate conjecture holds for $S$, it also 
  holds for $\tilde{S}$, cf. the introduction of \cite{mi}.
  Then the Tate conjecture also holds for $X$ since there is a finite morphism
  from a surface for which the Tate conjecture holds onto $X$, namely $\pi$.
  But the truth of the Tate conjecture implies the truth of the Artin--Tate
  conjecture by the main result of \cite{mi}.
  According to \cite{llr}, the results of \cite{mi} are also true in characteristic
  $p=2$.
  
  Since $\pi$ is purely inseparable, the zeta functions of $\tilde{S}$ and
  $X$ coincide. 
  In particular, we have $P_2(X,t)=P_2(S,t)$.
  Together with the fact that $\Br(\tilde{S})$ is trivial and
  $\alpha(\tilde{S})=g(C)\cdot g(F)$, we obtain the stated formula
  for ${\rm disc}\,\NS(X)$.
\qed

\section{Examples}
 
\subsection*{Bogomolov--Miyaoka--Yau}
For a surface over the complex numbers, this 
inequality states $c_1^2\leq9\chi$.
It is known to fail in positive characteristic, see e.g.
\cite[Section 3.4.1]{sz} or \cite[Kapitel 3.4.J]{bhh}.
Over the complex numbers, $\chi=1$
is the lowest value possible for a surface of general
type.
These surfaces have been studied for quite some time
and so it is interesting to note that also among them
there are counter-examples to this inequality in 
positive characteristic:

\begin{Theorem}
  \label{bmy}
  There exist surfaces of general type 
  with $\chi=1$ and $c_1^2=14$
  in characteristic $2$.
\end{Theorem}

\prf
We let $f:C\to\PP^1$ be an Artin--Schreier curve 
of genus $3$ as in Section \ref{curvesection}.
As rational vector field $\delta_C$ we choose the additive
vector field $\partial/\partial x$ from (\ref{singartin}),
which has a pole of order $4$ at infinity and no zeros.
On $F:=\PP^1$ we choose the vector field $\delta_F:=\delta_1$
from (\ref{vectorone}), which has a pole of order $4$ 
and three zeros of order $2$.

The quotient $(C\times F)/(\delta_C+\delta_F)$ has exactly
one singularity, which is elliptic of type $(19)_0$.
The resolution of singularities yields a surface of general
type with $\chi=1$ and $c_1^2=14$.
\qed


\begin{Remark}
  By Noether's formula $12\chi=c_1^2+c_2$, 
  this surface has negative $c_2$.
  Using an Artin--Schreier curve of genus $4$ curve 
  instead of genus $3$ in the previous construction,
  we obtain a surface of general type with $\chi=3$ and
  $c_1^2=30$, i.e. a counter-example to the
  Bogomolov--Miyaoka--Yau inequality with 
  small $\chi$ and $c_2$ positive.
\end{Remark}

\subsection*{Bounds on $\boldsymbol{(-2)}$-curves}
Over the complex numbers, a theorem of Miyaoka bounds
the number of disjoint $(-2)$-curves 
on a minimal surface of general type above
by $\frac{1}{9}(3c_2-c_1^2)$, 
cf. \cite[Section VII.4]{bhpv}.
This may fail in in positive characteristic
and Shepherd-Barron \cite[Theorem 4.1]{sbfol}
has shown that if there exist more than
$c_1^2+\frac{1}{2}c_2$ disjoint $(-2)$-curves then the surface
in question is uniruled.

However, there is usually a gap between these two bounds and 
the following theorem shows us  that there exist uniruled 
as well as non-uniruled surfaces in this gap.

\begin{Theorem}
  \label{minustwo}
  There exist minimal surfaces of general type in characteristic $2$ that 
  violate Miyaoka's bound on $(-2)$-curves that do not reach 
  Shepherd-Barron's bound.
  There exist uniruled as well as non-uniruled
  such surfaces.
\end{Theorem}

The non-uniruled surface that we present is birationally dominated by
an Abelian surface.\medskip

\prf
First, we give a non-uniruled example.
Let $E:=E_\alpha$ be an elliptic curve as in
Section \ref{curvesection} with $\alpha\neq0$, i.e.
$E$ is not supersingular.
The vector field $\delta_E:=\delta_{\alpha,1,\alpha}$ defined
by (\ref{singelliptic}) is an additive rational vector field on
$E$ with a zero and a pole of order $2$. 

We apply our construction to the data 
$(E,E,\delta_E+\delta_E)$.
The singular quotient 
$(E\times E)/(\delta_E+\delta_E)$ 
has two singularities of type $D_4$.

Resolving the singularities we obtain a minimal 
surface $X$ of general type with 
$\chi=1$, $c_1^2=4$ and hence $c_2=8$.
There are $6$ isolated $(-2)$-curves coming from the 
two $D_4$-singularities and hence Miyaoka's bound is
violated.

By construction, this surface is inseparably dominated 
by a blow-up of the Abelian surface $E\times E$.
Factoring the Frobenius morphism and blowing down
the $(-1)$-curves, we obtain a surjective 
morphism from $X$ onto the Abelian surface 
$(E\times E)^{(-1)}$.
If $X$ were uniruled, we would have a dominant
map $Y\to X$ from some birationally ruled surface $Y$.
Thus, there would exist a dominant map from $Y$ to
an Abelian surface which would have to factor over
the Albanese map of $Y$.
This is absurd since the image of the Albanese map of $Y$
is a curve.
Hence, $X$ is not uniruled.
\medskip

To obtain a uniruled example, we do the computations with the
example of Theorem \ref{bmy} but with an Artin--Schreier
curve of genus $2$ instead of genus $3$.
The singular quotient has exactly one singularity, which is
of type $D_8$.
The resolution is a surface $X$ with 
$\chi=1$, $c_1^2=8$ and $c_2=4$.
There are $5$ disjoint $(-2)$-curves coming from the resolution
of the $D_8$-singularity which is already enough to violate 
Miyaoka's bound.
On the other hand, this surface has $b_2=10$ and so the
rank of its N\'eron--Severi group is at most $10$
by the Igusa--Severi inequality.
Disjoint $(-2)$-curves are linearly independent in the
N\'eron--Severi group and they span a negative definite
lattice.
By the Hodge index theorem, there can be at most $9$ 
disjoint $(-2)$-curves on $X$.
Hence Shepherd-Barron's bound is not reached.
\qed

\section{Unbounded pathological behaviour}

\subsection*{Nonreduced Picard schemes}
In characteristic zero, Cartier's theorem states
that group schemes over a field are smooth.
This may fail in positive characteristic, 
and so one has to distinguish between the Picard scheme
and the Picard variety of a given variety.
A first example of a smooth variety with non-reduced Picard
scheme was found by Igusa \cite{ig2}.
Whereas $\frac{1}{2}b_1$ gives the dimension of
the Picard variety, the number $h^{01}$ gives
the dimension of the tangent space to the Picard
scheme and hence we only have
an inequality $\frac{1}{2}b_1 \leq h^{01}$, with
equality if and only if the Picard scheme is reduced.

\begin{Theorem}
  \label{picard}
  Given an integer $q\geq2$, there exists a family 
  $\{X_i\}_{i\in\NN}$ of
  uniruled surfaces of general type in
  characteristic $2$
  all having the same Picard variety 
  of dimension $q$ such that
  $$
   h^{01}(X_i)\,=\,h^1(\OO_{X_i})\,\to\,\infty
   \,\mbox{ as }\,i\to\infty\,.
  $$
  Thus, the Picard scheme can get arbitrarily non-reduced,
  even when fixing the Picard variety.
\end{Theorem}

\prf
We let $h:=q+1$ and let $\varphi:C\to\PP^1$ be the 
Artin--Schreier cover of $\PP^1$
given by $z^2-z=x^{2h-1}$, which defines a curve
of genus $q$ as explained in Section \ref{curvesection}.
Let $m$ be the largest integer less or equal to
$\frac{1}{2}(q-1)$.
We choose $m$ distinct elements
$\{a_i\}_{i=1,...,m}$ of the ground
field $k$ and consider the rational vector field
$$
\delta_C \,:=\, \varphi^{\ast} \left(
\prod_{i=1}^{m} \frac{1}{(x-a_i)^2}\,\cdot\,
\frac{\partial}{\partial x}
\right)\,.
$$
This rational vector field is additive and has $(g-1)$ poles of
order $2$ on $C$.

We let $F:=\PP^1$ and choose additive vector fields
with poles and zeros of order $2$.
Clearly, we can find a family $\{\delta_F^i\}_{i\in\NN}$ 
of such vector fields such that the degree
$d_F^i$ of its divisor of poles tends to infinity.
For example, we can use vector fields as in
(\ref{singrational}).

We set $\delta^i:=\delta_C+\delta_F^i$ on $S:=C\times F$
and let $X_i$ be the uniruled surface constructed from data
$(C,F,\delta^i)$. 
By Proposition \ref{bettiinvariants}, the Albanese variety, i.e.
the dual of the Picard variety,
of the $X_i$'s is the Jacobian of $C^{(-1)}$.

Since $\delta_C$ and $\delta_F^i$ have only poles of order $2$,
the singular quotient $X_i':=S/\delta^i$ is a normal surface,
which has only Du~Val singularities of type $D_4$.
In particular, we can compute $h^{01}(X_i)$ on $X_i'$.
By Proposition \ref{hodgeinvariants}, this number equals
$g(C)+h^1({\cal N}_i^\vee)$,
where ${\cal N}_i$ is given by
$$
{\cal N}_i^{\otimes (-2)} \,\iso\, \OO_F(-d_F^i-2)\,\boxtimes\,
\OO_C
$$
In Section \ref{curvesection} we noted
that the $2$-rank of the Artin--Schreier curve $C$ is zero, 
which implies that
$$
{\cal N}_i^\vee \,\iso\, \OO_F(-\frac{d_F^i}{2}-1)\,\boxtimes\,
\OO_C\,.
$$
Using the K\"unneth formula, we see that
$h^1({\cal N}_i^\vee)$ tends to infinity as 
$d_F^i$ tends to infinity.
\qed

\subsection*{Global $\mathbf1$-forms}
In characteristic zero, Hodge theory implies that every
global $1$-form on a variety is the pull-back of a $1$-form
from its Albanese variety via the Albanese map.
Igusa \cite{ig} showed that the pull-back of a non-trivial global
$1$-form from the Albanese variety of a smooth variety via the Albanese
map remains non-trivial in arbitrary characteristic.
Therefore, one always has the inequality 
$h^{10}(X)=h^0(\Omega_X^1)\geq h^0(\Omega_{\Alb(X)}^1)=\frac{1}{2}b_1(X)$.
On the other hand, Igusa \cite{ig2} also gave an example of a
surface in positive characteristic with strict inequality.

\begin{Theorem}
  \label{albanese}
  Given an integer $q\geq2$, there exists a family 
  $\{X_i\}_{i\in\NN}$ of
  uniruled surfaces of general type in
  characteristic $2$
  all having the same Albanese variety 
  of dimension $q$ such that
  $$
   h^{10}(X_i)\,=\,h^0(\Omega_{X_i}^1)\,\to\,\infty
   \,\mbox{ as }\,i\to\infty\,.
  $$
\end{Theorem}

\prf
We take the family of surfaces from Theorem \ref{picard}.
The computations from the proof and the long exact sequence
of cohomology applied to (\ref{dclosedforms}) yield
$h^{10}(X_i)\geq g(C)+3d_F^i/2-1$, 
which tends to infinity as $i$ tends to infinity.
\qed

\subsection*{Hodge symmetries}
In characteristic zero one has not only that the Fr\"olicher
spectral sequence of a projective variety always degenerates at 
$E_1$-level but also the Hodge symmetries $h^{ij}=h^{ji}$.
In positive characteristic, these symmetries can be violated
even if the Fr\"olicher spectral sequence degenerates
at $E_1$-level, cf. \cite[Proposition 16]{semex}
and \cite[Remarque 2.6 (ii)]{delill}.

\begin{Theorem}
  \label{symmetry}
  Given an integer $q\geq2$, there exists a family 
  $\{X_i\}_{i\in\NN}$ of
  uniruled surfaces of general type in
  characteristic $2$
  all having the same Albanese variety 
  of dimension $q$ such that
   $$
   h^{10}(X_i)-h^{01}(X_i)\,\to\,\infty
   \,\mbox{ as }\,i\to\infty\,,
  $$
  i.e. the Hodge symmetries fail.
\end{Theorem}

\prf
We use the family used in the proofs of Theorem \ref{picard} 
and Theorem \ref{albanese}, where we have seen
$h^{01}(X_i)=g(C)+d_F^i/2$ and
$h^{10}(X_i)\geq g(C)+d_F^i-1$.
\qed

\subsection*{Global vector fields}
Over the complex numbers, a surface of general type has no
non-trivial global vector fields.
This follows from the fact that the vector space of 
global vector fields can be identified with the tangent 
space to the group of biholomorphic automorphisms, 
and this group is finite for a surface of general type.
Counter-examples in positive characteristic can 
be found, e.g.
in \cite{lang} or in \cite{sbfol}, where first examples
of non-uniruled surfaces of general type with vector fields
appear.

\begin{Proposition}
  \label{globalvectorfields}
  Let $X$ be a surface constructed from data $(C,F,\delta)$
  in characteristic $2$.
  We assume that $g(C)\leq1$, and $g(F)\leq1$ as well as $p_g(X)\neq0$.
  Then $X$ possesses non-trivial global vector fields.
\end{Proposition}

\prf
We keep the notations of the proof of Proposition \ref{hodgeinvariants}.
We consider the morphism $\varpi:X\to\tilde{S}^{(-1)}$.
Dualising (\ref{cotangent}) and plugging in (\ref{relative})
we obtain an exact sequence
\begin{equation}
  \label{vector}
0\,\to\,\omega_X\otimes\varpi^\ast\omega_{\tilde{S}^{(-1)}}^\vee\,\to\,
\Theta_X^1\,\to\,\varpi^{\ast}\Theta_{\tilde{S}^{(-1)}}^1\,\to\,...
\end{equation}
The assumptions on $g(C)$, $g(F)$ and $p_g(X)$ make sure that
$\omega_X$ and $\omega_{\tilde{S}^{(-1)}}^\vee$ have non-trivial global
sections.
Hence 
$\omega_X\otimes\varpi^\ast\omega_{\tilde{S}^{(-1)}}^\vee$ 
has non-trivial global sections, 
which yields non-trivial global vector fields via (\ref{vector}).
\qed

\begin{Theorem}
  \label{unboundedvectorfields}
  In characteristic $2$, there exist families of surfaces of general
  type where $c_1^2$ and $\chi$ tend to infinity and such that
  each member
  of this family possesses non-trivial global vector fields.
  Moreover, we can find such families in which every member is
  uniruled, resp. not uniruled.
\end{Theorem}

\prf
First, we give unirational examples.
We set $F:=\PP^1$ and choose a family of additive rational
vector fields $\delta_F^i$ poles and zeros of order $2$ such
that $d_F^i$, the degree of the divisor of poles of $\delta_F^i$,
tends to infinity as $i$ tends to infinity.
For example, we could use vector fields of the form (\ref{singrational}).

The quotient $X_i':=(F\times F)/(\delta^i_F+\delta^i_F)$ is a normal
surface with Du~Val singularities of type $D_4$.
If $X_i$ is the surface constructed from $(F,F,\delta_F^i+\delta_F^i)$
then its invariants $\chi$ and $K^2$ coincide with those of $X_i'$.
Thus it is enough to show that $K_{X_i'}^2$ and $\chi(\OO_{X_i'})$
are unbounded as $i$ tends to infinity.
However, this can easily be seen from
Proposition \ref{canonicaldivisor} and
Proposition \ref{hodgeinvariants}.

By Proposition \ref{globalvectorfields}, the surfaces $X_i$
possess non-trivial global vector fields.\medskip

To construct a family of surfaces of general type that is
not uniruled, we let $\varphi:E\to\PP^1$ be an elliptic curve
given as a separable double cover branched over
$x=0$ and $x=\infty$.
We choose $2n$ pairwise distinct and non-zero elements
$\{a_i, b_i\}_{i=1,...,n}$ of the ground
field $k$ and consider
$$
\delta_E^n \,:=\, \varphi^{\ast} \left(
\prod_{i=1}^{n} \frac{(x-a_i)^2}{(x-b_i)^2}\,\cdot\,
\frac{\partial}{\partial x}
\right)\,.
$$
This defines an additive rational vector field on $E$ with
poles and zeros of order $2$.

Arguing as before we see that
the surfaces constructed from data $(E,E,\delta_E^n+\delta_E^n)$
have non-trivial global vector fields and
that $\chi$ and $c_1^2$ tend to infinity as $n$ tends to
infinity.
Arguing as in the proof of Theorem \ref{minustwo} we see
that these surfaces are not uniruled.
\qed


\subsection*{Inseparability of the canonical map}
Curves of general type, i.e. of genus at least $2$, fall into
two classes: the hyperelliptic and the non-hyperelliptic ones.
Whether a curve is hyperelliptic or not can be seen from
the canonical map: it is either a separable morphism of 
degree $2$ onto $\PP^1$ or defines an embedding.
In any case, it is a finite and separable morphism
onto its image.

For surfaces, the canonical map may be empty and is usually
only a rational map.
The canonical map of a product of two curves of general
type is a finite morphism onto a rational,
a ruled or a general type surface, depending on whether
the curves are hyperelliptic or not.
In these examples, the canonical map is a separable
morphism.
We will see in the next theorem that the canonical
map can become inseparable and that this is not a
sporadic phenomenon.

It follows from Shepherd-Barrons work
\cite[Theorem 27]{sb} that for $c_1^2$ and $\chi$ sufficiently
large, $|3K_X|$ defines a birational morphism.
Also, if $X$ has no pencil of curves of arithmetic
genus $2$ then already $|2K_X|$ defines a birational morphism
if $c_1^2$ and $\chi$ are sufficiently large.
However, the surfaces presented in Theorem \ref{canonicalinsep}
do not possess pencils of curves of small arithmetic genus.
Hence the inseparability of the canonical map is not related
to the existence of special fibrations of low genus.

\begin{Theorem}
  \label{canonicalinsep}
  In characteristic $2$, there exist families of surfaces of general
  type where $c_1^2$ and $\chi$ tend to infinity and such that
  the canonical map of each member is a generically
  finite and inseparable morphism onto a rational surface.
  
  Moreover, given a natural number $b$, we can find such families
  that do not possess pencils of curves of arithmetic genus
  less than $b$.
\end{Theorem}

\prf
We consider the family of unirational surfaces of general type
constructed in Theorem \ref{unboundedvectorfields} and use
the notations from the proof.

We will work on the singular surfaces $X_i'$, which have only
Du~Val singularities of type $D_4$.
The $X_i$'s are $\alpha_{{\cal L}_i}$-torsors 
over $S^{(-1)}:=(\PP^1\times\PP^1)^{(-1)}$.
It is not difficult to see that ${\cal L}_i$ and 
$\omega_{S^{(-1)}}\otimes{\cal L}_i$ are very ample
line bundles on $S^{(-1)}$.
From (\ref{dualisingsheaf}) it follows that the canonical
sheaf on $X_i'$ is ample.
Hence, $X_i'$ is the canonical model of $X_i$ and the canonical
map factors over $X_i'$ which justifies to work with the
$X_i'$'s rather than the $X_i$'s.

Let us analyse the canonical system of the $X_i$'s.
We have already seen above that there is a finite and 
purely inseparable morphism $\varpi_i:X_i'\to S^{(-1)}$
which has the structure of an $\alpha_{ {\cal L}_i}$-torsor.
Pushing $\omega_{X_i'}$ forward to
$S^{(-1)}$ and using (\ref{dualisingsheaf}),
we obtain an extension
$$
0\,\to\,\omega_{S^{(-1)}}\otimes{\cal L}_i\,\to\,
\varpi_{i,\ast}\, (\omega_{X_i'}) \,\to\,
\omega_{S^{(-1)}}
\,\to\,0 \,.
$$
The line bundle $\omega_{S^{(-1)}}$ has no global sections
and it is not difficult to see that the line bundle
$\omega_{S^{(-1)}}\otimes {\cal L}_i$ defines an embedding
$\varphi_i$ of $S^{(-1)}$ into $\PP^{p_g(X_i')-1}$.
Thus, the canonical map of $X_i'$ factors as
$\varphi_i\circ\varpi_i$.
In particular, it is a finite and
inseparable morphism onto a rational surface.

Suppose $X_i$ has a pencil of curves of arithmetic genus $d_i$, 
say with generic fibre $D_i$.
Then the adjunction formula yields $2d-2=K_{X_i} D_i$ since a
fibre has $D_i^2=0$.
The canonical divisor on $X_i$ is the pull-back of a divisor
on $\tilde{S}^{(-1)}$ by (\ref{dualisingsheaf}).
Using the projection formula and (\ref{linebundle}), 
it is not difficult to see that if $d_F^i$ is sufficiently large 
then $K_{X_i} D_i > 2b-2$ for any given bound $b$.
Hence, if $d_F^i$ is sufficiently large, then $X_i$ does not
possess a pencil of curves of arithmetic genus less than
$b$.
\qed

\end{document}